# Invariance principle, multifractional Gaussian processes and long-range dependence

Serge Cohen[a] and Renaud Marty[b]

[a]*Laboratoire de Statistique et Probabilités, Université Paul Sabatier, 118, Route de Narbonne, 31062 Toulouse Cedex 4, France. E-mail: Serge.Cohen@math.ups-tlse.fr*

[b]*Department of Mathematics, 103 Multipurpose Science and Technology Building, University of California at Irvine, Irvine, CA 92697-3875, USA. E-mail: renaud.marty@iecn.u-nancy.fr*

**Abstract.** This paper is devoted to establish an invariance principle where the limit process is a multifractional Gaussian process with a multifractional function which takes its values in $(1/2, 1)$. Some properties, such as regularity and local self-similarity of this process are studied. Moreover the limit process is compared to the multifractional Brownian motion.

**Résumé.** Ce papier a pour but d'établir un principe d'invariance dont le processus limite est gaussien et multifractionnaire avec une fonction de Hurst à valeurs dans $(1/2, 1)$. Des propriétés telles que la régularité et l'autosimilarité locale de ce processus sont étudiées. De plus, le processus limite est comparé au mouvement brownien multifractionnaire.

*AMS 2000 subject classifications.* 60F17; 60G15

*Keywords:* Invariance principle; Long range dependence; Multifractional process; Gaussian processes

## 1. Introduction

Fractional and multifractional processes or fields have been extensively studied because they provide relevant models in many situations such as mathematical finance, network traffic, physics, and other fields related to applied mathematics. See, e.g., [10] for a convenient reference.

The most famous and simplest fractional process is fractional Brownian motion which has been introduced by Mandelbrot and Van Ness [12]. For a fixed $H \in (0,1)$, it can be defined as a centered Gaussian process $W_H$ with covariance

$$\mathbb{E}[W_H(t)W_H(s)] = \frac{\mathbb{E}[W_H(1)^2]}{2}(|t|^{2H} + |s|^{2H} - |t-s|^{2H}).$$

Its main properties, which depend on $H$, are the following.

- $W_H$ is self-similar with index $H$.

- $W_H$ admits $H$ as Hölder exponent at each point.
- If $H > 1/2$ the increments of $W_H$ (which are stationary) satisfy the long-range dependence relation as $n \to \infty$

$$\mathbb{E}[W_H(1)(W_H(n+1) - W_H(n))] \sim H(2H-1)\mathbb{E}[W_H(1)^2]n^{2H-2}. \tag{1}$$

Let us focus on this last property. Notice that the sequence of the increments of fractional Brownian motion called the fractional white noise is not the only Gaussian sequence with long-range dependence (see, for instance, [15], p. 336). However, the fractional white noise serves as a universal Gaussian model for long-range stationary phenomena. This is due to the invariance principle, which was established in [9, 15].

**Theorem 1 (Invariance principle).** *Let $H$ belong to $(1/2, 1)$, and $\{X_n\}_{n \in \mathbb{N}}$ be a stationary sequence of centered Gaussian random variables with covariance satisfying when $n \to \infty$*

$$\mathbb{E}[X_0 X_n] \sim c n^{2H-2},$$

*with $c > 0$. Then defining for every $t > 0$,*

$$S^N(t) = \frac{1}{N^H} \sum_{n=1}^{\lfloor Nt \rfloor} X_n,$$

*the finite-dimensional distributions of $S^N$ converge to those of $c_0 W_H$ as $N$ goes to $\infty$, where $W_H$ is a fractional Brownian motion of index $H$ and $c_0^2 = H^{-1}(2H-1)^{-1}c$.*

In the previous statement, $\lfloor x \rfloor$ denotes the integer part of the real number $x$.

One of the main drawbacks of fractional Brownian motion for applications lies in the homogeneity of its properties, due to the fact that its pointwise Hölder exponent is constant. Hence, multifractional processes and fields, for instance, the multifractional Brownian motion, have been introduced and attracted attention [3, 4, 6, 8, 11, 13]. Multifractional processes have locally the same properties as the fractional processes, but not globally. In fact their properties are not governed by a constant exponent $H \in (0, 1)$, but by a $(0, 1)$-valued function $h$ which is called the multifractional function. For instance, multifractional processes are locally self-similar and their pointwise Hölder exponents vary along the trajectory.

Because they take into account the variations of properties such as regularity, multifractional processes with a $(1/2, 1)$-valued multifractional function could be relevant alternatives to fractional processes with $H \in (1/2, 1)$ to provide models for long-range phenomena [2]. The main aim of this paper is to prove the existence of multifractional Gaussian processes which can serve as universal Gaussian models for long-range dependence. More precisely, we establish an invariance principle where the limit process is a multifractional Gaussian process with long-range increments. A study of these processes is developed. Besides local self-similarity and Hölder regularity, a representation result is given.

Maybe the most famous multifractional process is the multifractional Brownian motion. In this paper, we also compare as much as possible the limit processes we obtain with the multifractional Brownian motion.

In Section 2, we recall the definition and some properties of the multifractional Brownian motion. The main result of this paper is established in Section 3. Section 4 deals with local self-similarity and regularity properties of the limit processes which are obtained in the main result. In Section 5, we give a representation of the limit processes. This representation allows us to justify a universal property for the limit processes. Finally, we give in Section 6 some applications of our results, in particular, an extension of results of [14]. The proof of some technical lemmas are postponed to the Appendices.

*Notation*

- For two random processes (or fields) $Z_1$ and $Z_2$, we denote by $Z_1 =^{\mathcal{D}} Z_2$, the fact that the finite dimensional margins of $Z_1$ are equal in distribution to those of $Z_2$.

- For a random process (or field) $\widetilde{Z}$ and a sequence $\{Z_N\}_{N\in\mathbb{N}}$ (resp. $\{Z_\varepsilon\}_{\varepsilon>0}$) of random processes (or fields), we denote by $\lim^{\mathcal{D}}_{N\to\infty} Z_N = \widetilde{Z}$ (resp. $\lim^{\mathcal{D}}_{\varepsilon\to 0} Z_\varepsilon = \widetilde{Z}$) the fact that the finite-dimensional distributions of $Z_N$ (resp. $Z_\varepsilon$) converge to those of $\widetilde{Z}$ as $N$ goes to $\infty$ (resp. $\varepsilon$ goes to 0).

## 2. Some preliminaries on multifractional Gaussian processes

Multifractional Brownian motion was the first multifractional Gaussian process introduced, independently in [6] and [13].

Let $B(\mathrm{d}\xi)$ be a Brownian measure and $\widehat{B}(\mathrm{d}x)$ be its Fourier transform (See [15], Chapter 7 for more details). Multifractional Brownian motion with a multifractional function $h\colon [0,\infty) \to (0,1)$ can be defined for every $t$ by

$$W_h(t) = \frac{1}{C(h(t))} \int_{-\infty}^{\infty} \frac{\mathrm{e}^{-\mathrm{i}xt} - 1}{|x|^{h(t)+1/2}} \widehat{B}(\mathrm{d}x) \tag{2}$$

following [6]. The constant $C(h(t))$ is such that $\mathbb{E}[W_h(1)^2] = 1$. In this case, the covariance of the multifractional Brownian motion is given for every $t$ and $s$ by [2]

$$\mathbb{E}[W_h(t)W_h(s)] = D(h(t),h(s))(|t|^{h(t)+h(s)} + |s|^{h(t)+h(s)} - |t-s|^{h(t)+h(s)}), \tag{3}$$

where for every $H_1, H_2$

$$D(H_1, H_2) = \frac{\sqrt{\Gamma(2H_1+1)\Gamma(2H_2+1)\sin(\pi H_1)\sin(\pi H_2)}}{2\Gamma(H_1+H_2+1)\sin(\pi(H_1+H_2)/2)}.$$

Now we assume that $h$ is $\beta$-Hölder continuous and

$$\sup h < \beta. \tag{4}$$

Then multifractional Brownian motion has locally the same properties as fractional Brownian motion. Multifractional Brownian motion is locally self-similar, that is, for every $t \geq 0$,

$$\lim_{\varepsilon \to 0} \left\{ \frac{W_h(t+\varepsilon u) - W_h(t)}{\varepsilon^{h(t)}} \right\}_{u \geq 0} = \{W_{h(t)}(u)\}_{u \geq 0},$$

where the convergence is in distribution in the space of continuous functions endowed with the topology of the uniform convergence on compact sets, and for every $t$ its Hölder pointwise exponent $\alpha_{W_h}(t)$ is almost surely equal to $h(t)$.

## 3. The main result

In this paper, we consider a centered Gaussian field $X = \{X_n(H), H \in (1/2, 1), n \in \mathbb{N}\}$, and two real numbers $a$ and $b$ such that $1/2 < a \leq b < 1$. We also consider a continuous function $h\colon \mathbb{R} \to [a,b] \subset (1/2, 1)$. For every $n, N \in \mathbb{N}$ and $t > 0$, we define $h_n^N = h(n/N)$ and

$$S_h^N(t) = \sum_{n=1}^{\lfloor Nt \rfloor} \frac{X_n(h_n^N)}{N^{h_n^N}}.$$

Our aim is to study the asymptotic behavior of $\{S_h^N\}_{N\in\mathbb{N}}$ as $N$ goes to $\infty$. We will use the following assumptions.

- **Assumption (i)** For every $M > 0$, the map

$$(j, k, H_1, H_2) \mapsto \mathbb{E}[X_j(H_1) X_k(H_2)]$$

is bounded on $\{(j,k) \in \mathbb{N}^2, |j - k| \leq M\} \times [a, b]^2$.
- **Assumption (ii)** There exists a continuous function $R : [a, b]^2 \to (0, \infty)$ such that

$$\lim_{j-k \to \infty} \sup_{(H_1, H_2) \in [a,b]^2} |(j-k)^{2-H_1-H_2} \mathbb{E}[X_j(H_1) X_k(H_2)] - R(H_1, H_2)| = 0. \tag{5}$$

Assumption (ii) expresses that the field is asymptotically stationary and fractional for each $H$. Hence, for a fixed $H$, the random sequence $n \mapsto X_n(H)$ satisfies the classical invariance principle.

For a centered Gaussian field, $X = \{X_n(H), H \in (1/2, 1), n \in \mathbb{N}\}$ satisfying Assumption (ii); the function $R$ in Assumption (ii) will be called the asymptotic covariance of $X$. More generally, a measurable and bounded function $R : [a, b]^2 \to (0, \infty)$ such that there exists a centered Gaussian field $X = \{X_n(H), H \in (1/2, 1), n \in \mathbb{N}\}$ satisfying Assumption (ii) with $R$ as asymptotic covariance will be called an asymptotic covariance.

Fix an asymptotic covariance $R$. We denote by $\mathcal{G}_R$ the space of all centered Gaussian fields $\{X_n(H)\}_{n,H}$ satisfying Assumptions (i) and (ii) with $R$ as asymptotic covariance. We define for every $t, s, H_1, H_2$

$$\mathcal{R}(t, s; H_1, H_2) = R(H_1, H_2) 1_{t \geq s} + R(H_2, H_1) 1_{t < s}, \tag{6}$$

and the function

$$\mathcal{R}^* : [0, \infty)^2 \to \mathbb{R},$$
$$(\theta, \sigma) \mapsto \mathcal{R}(\theta, \sigma; h(\theta), h(\sigma)) |\theta - \sigma|^{h(\theta) + h(\sigma) - 2}. \tag{7}$$

Let us remark the symmetric property fulfilled by $\mathcal{R}$: for every $t, s, H_1, H_2$,

$$\mathcal{R}(t, s; H_1, H_2) = \mathcal{R}(s, t; H_2, H_1).$$

Let us also note that $\mathcal{R}^*$ is locally integrable. Indeed, by the boundedness of $R$, on every compact set $K \subset [0, \infty)^2$ there exists a constant $c = c_K$, such that for every $(\theta, \sigma) \in K$ satisfying $\theta \neq \sigma$

$$|\mathcal{R}^*(\theta, \sigma)| \leq c \times (|\theta - \sigma|^{2a-2} + |\theta - \sigma|^{2b-2}). \tag{8}$$

**Theorem 2.** *Under Assumptions* (i) *and* (ii), *when $N$ goes to $\infty$, the finite-dimensional distributions of $S_h^N$ converge to those of a centered Gaussian process $\widetilde{S}_h$ with covariance given for $t, s \geq 0$ by:*

$$\mathbb{E}[\widetilde{S}_h(t) \widetilde{S}_h(s)] = \int_0^t d\theta \int_0^s d\sigma \, \mathcal{R}(\theta, \sigma; h(\theta), h(\sigma)) |\theta - \sigma|^{h(\theta) + h(\sigma) - 2}$$
$$= \int_0^t d\theta \int_0^s d\sigma \, \mathcal{R}^*(\theta, \sigma), \tag{9}$$

*where the integral in the right-hand side of* (9) *is always defined.*

**Proof.** For every $n$ and $N$, we define $X_{n,N} := X_n(h_n^N)/N^{h_n^N}$ and we let for every $t, s \geq 0$

$$I(s, t) = \int_0^t d\theta \int_0^s d\sigma \, \mathcal{R}^*(\theta, \sigma). \tag{10}$$

Because $\{X_n(H)\}_{n,H}$ is a centered Gaussian field, it is enough to show that for every $t, s \geq 0$:

$$\lim_{N \to \infty} \mathbb{E}[S_h^N(t) S_h^N(s)] = \lim_{N \to \infty} \sum_{j=1}^{\lfloor Nt \rfloor} \sum_{k=1}^{\lfloor Ns \rfloor} \mathbb{E}[X_{k,N} X_{j,N}] = I(t, s).$$

Formula (5) deals with the asymptotic behavior of $\mathbb{E}[X_j(H_1)X_k(H_2)]$ when $j-k$ goes to $+\infty$. Moreover, it also gives

$$\lim_{j-k\to-\infty} \sup_{(H_1,H_2)\in[a,b]^2} |(k-j)^{2-H_1-H_2}\mathbb{E}[X_j(H_1)X_k(H_2)] - R(H_2,H_1)| = 0, \tag{11}$$

so that we can write for every $j$, $k$ (with $j \neq k$), $H_1$ and $H_2$

$$\mathbb{E}[X_j(H_1)X_k(H_2)] = |j-k|^{H_1+H_2-2}(\mathcal{R}(j,k;H_1,H_2) + r_0(j,k;H_1,H_2)), \tag{12}$$

where

$$\lim_{|j-k|\to\infty} \sup_{(H_1,H_2)\in[a,b]^2} r_0(j,k;H_1,H_2) = 0.$$

Let $\eta > 0$. Following (12) and the fact that $a \leq h \leq b$, there exists an integer $M = M_\eta$ such that for $|j-k| > M$ and for every $N$,

$$\frac{1}{N^2}\mathcal{R}^*\left(\frac{j}{N},\frac{k}{N}\right) - \frac{\eta}{N^2}\left|\frac{j-k}{N}\right|^{h_j^N+h_k^N-2}$$

$$\leq \mathbb{E}[X_{k,N}X_{j,N}] \leq \frac{1}{N^2}\mathcal{R}^*\left(\frac{j}{N},\frac{k}{N}\right) + \frac{\eta}{N^2}\left|\frac{j-k}{N}\right|^{h_j^N+h_k^N-2}. \tag{13}$$

Thanks to Lemma 2 (see the Appendix),

$$I(s,t) - \eta J(s,t) \leq \liminf_{N\to\infty} \sum_{j=1}^{\lfloor Nt \rfloor}\sum_{k=1}^{\lfloor Ns \rfloor} \mathbb{E}[X_{k,N}X_{j,N}]1_{|j-k|>M}$$

$$\leq \limsup_{N\to\infty} \sum_{j=1}^{\lfloor Nt \rfloor}\sum_{k=1}^{\lfloor Ns \rfloor} \mathbb{E}[X_{k,N}X_{j,N}]1_{|j-k|>M} \leq I(s,t) + \eta J(s,t), \tag{14}$$

where

$$J(s,t) = \int_0^t d\theta \int_0^s d\sigma\, |\theta-\sigma|^{h(\theta)+h(\sigma)-2}. \tag{15}$$

Now we deal with $\sum_{j=1}^{\lfloor Nt \rfloor}\sum_{k=1}^{\lfloor Ns \rfloor} \mathbb{E}[X_{k,N}X_{j,N}]1_{|j-k|\leq M}$. From Assumption (i) and $a \leq h \leq b$, we get for every $j$, $k$ and $N$

$$|\mathbb{E}[X_{k,N}X_{j,N}]| \leq \frac{1}{N^{h_j^N+h_k^N}} \sup_{(H_1,H_2)\in[a,b]^2} |\mathbb{E}[X_j(H_1)X_k(H_2)]|$$

$$\leq \frac{1}{N^{2a}} \sup_{(H_1,H_2)\in[a,b]^2} |\mathbb{E}[X_j(H_1)X_k(H_2)]|, \tag{16}$$

where $\sup_{(H_1,H_2)\in[a,b]^2} |\mathbb{E}[X_j(H_1)X_k(H_2)]| < \infty$. Moreover,

$$\sum_{j=1}^{\lfloor Nt \rfloor}\sum_{k=1}^{\lfloor Ns \rfloor} 1_{|j-k|\leq M} \leq [N(t+s)](1+2M). \tag{17}$$

Hence, combining (16) with (17) and using $2a > 1$, we get

$$\lim_{N\to\infty} \sum_{j=1}^{\lfloor Nt \rfloor}\sum_{k=1}^{\lfloor Ns \rfloor} \mathbb{E}[X_{k,N}X_{j,N}]1_{|j-k|\leq M} = 0. \tag{18}$$

Combining with (14) and letting $\eta \to 0$, we obtain for every $t$ and $s$,

$$\lim_{N \to \infty} \mathbb{E}[S_h^N(t) S_h^N(s)] = I(t,s).$$

Hence, the finite-dimensional distributions of $S_h^N$ converge to those of a centered Gaussian process $\widetilde{S}_h$ with covariance $\mathbb{E}[\widetilde{S}_h(t)\widetilde{S}_h(s)] = I(t,s)$. $\square$

## 4. Properties of the limit process

In this section, we study some properties of the limit process obtained in Theorem 2. First we deal with local self-similarity.

**Proposition 1.** *We assume that $h$ is Hölder continuous. Then the process $\widetilde{S}_h$ is locally self-similar, more precisely*

$$\lim_{\varepsilon \to 0}^{\mathcal{D}} \left\{ \frac{\widetilde{S}_h(t+\varepsilon u) - \widetilde{S}_h(t)}{\varepsilon^{h(t)}} \right\}_{t,u \geq 0} = \{\mathcal{T}(t,u)\}_{t,u \geq 0} = \mathcal{T}, \qquad (19)$$

*and the tangent process $\mathcal{T}$ is the centered Gaussian field with covariance*

$$\mathbb{E}[\mathcal{T}(t,u)\mathcal{T}(s,v)] = \begin{cases} \frac{R(h(t),h(t))}{4h(t)^2 - 2h(t)}(|u|^{2h(t)} + |v|^{2h(t)} - |u-v|^{2h(t)}) & \text{if } t = s, \\ 0 & \text{if } t \neq s. \end{cases}$$

*Moreover, for every $t$, as $\varepsilon \to 0$, the field $\{(\widetilde{S}_h(t+\varepsilon u) - \widetilde{S}_h(t))/\varepsilon^{h(t)}\}_{u \geq 0}$ converges in distribution in the space of continuous functions endowed with the topology of the uniform convergence to $\mathcal{T}(t,\cdot)$ which is a fractional Brownian motion, such that $\mathbb{E}[\mathcal{T}(t,1)^2] = R(h(t),h(t))/(2h(t)^2 - h(t))$ with Hurst index $h(t)$.*

**Proof.** We have

$$\mathbb{E}\left[\frac{(\widetilde{S}_h(t+\varepsilon u) - \widetilde{S}_h(t))(\widetilde{S}_h(s+\varepsilon v) - \widetilde{S}_h(s))}{\varepsilon^{h(t)+h(s)}}\right]$$

$$= \frac{1}{\varepsilon^{h(t)+h(s)}} \int_t^{t+\varepsilon u} d\theta \int_s^{s+\varepsilon v} d\sigma\, \mathcal{R}^*(\theta,\sigma)$$

$$= \frac{1}{\varepsilon^{h(t)+h(s)-2}} \int_0^u d\theta \int_0^v d\sigma\, \mathcal{R}^*(t+\varepsilon\theta, s+\varepsilon\sigma)$$

$$= \int_0^u d\theta \int_0^v d\sigma\, \mathcal{R}^{\#}(t,s,\theta,\sigma,\varepsilon),$$

where

$$\mathcal{R}^{\#}(t,s,\theta,\sigma,\varepsilon) = \frac{1}{\varepsilon^{h(t)+h(s)-2}} \mathcal{R}^*(t+\varepsilon\theta, s+\varepsilon\sigma).$$

First, we assume that $t = s$. We have

$$\mathcal{R}^{\#}(t,t,\theta,\sigma,\varepsilon) = \varepsilon^{h(t+\varepsilon\theta)+h(t+\varepsilon\sigma)-2h(t)}$$
$$\times R(t+\varepsilon\theta, t+\varepsilon\sigma, h(t+\varepsilon\theta), h(t+\varepsilon\sigma))|\theta - \sigma|^{h(t+\varepsilon\theta)+h(t+\varepsilon\sigma)-2}. \qquad (20)$$

Since $\mathcal{R}$ is continuous, when $\varepsilon \to 0$ we have for $\theta \neq \sigma$

$$R(t+\varepsilon\theta, t+\varepsilon\sigma, h(t+\varepsilon\theta), h(t+\varepsilon\sigma))|\theta - \sigma|^{h(t+\varepsilon\theta)+h(t+\varepsilon\sigma)-2}$$
$$\to R(h(t),h(t))|\theta - \sigma|^{2h(t)-2}. \qquad (21)$$

Moreover, $h$ is Hölder continuous, so there exist $\alpha_h > 0$ and $c_h > 0$ such that for every $\theta \in [0, u]$ and $\sigma \in [0, v]$

$$|h(t + \varepsilon\theta) - h(t)| \leq c_h \varepsilon^{\alpha_h} u^{\alpha_h} \quad \text{and} \quad |h(t + \varepsilon\sigma) - h(t)| \leq c_h \varepsilon^{\alpha_h} v^{\alpha_h}. \tag{22}$$

Hence, when $\varepsilon \to 0$ we have

$$\varepsilon^{h(t+\varepsilon\theta)+h(t+\varepsilon\sigma)-2h(t)} = \exp((h(t + \varepsilon\theta) + h(t + \varepsilon\sigma) - 2h(t))\log(\varepsilon)) \to 1. \tag{23}$$

Combining (20), (21) and (23) we obtain that for every $t$ and almost every $\theta$ and $\sigma$, when $\varepsilon \to 0$ we have

$$\mathcal{R}^{\#}(t, t, \theta, \sigma, \varepsilon) \to R(h(t), h(t))|\theta - \sigma|^{2h(t)-2}. \tag{24}$$

Note that because (22) the convergence (23) is uniform in $(\theta, \sigma) \in [0, u] \times [0, v]$. Then using (20) and the fact that $R$ is bounded, we prove that there exists a constant $c$ such that for every $\varepsilon$, $t$ and almost every $(\theta, \sigma) \in [0, u] \times [0, v]$ we have

$$\mathcal{R}^{\#}(t, t, \theta, \sigma, \varepsilon) \leq c \times (|\theta - \sigma|^{2b-2} + |\theta - \sigma|^{2a-2}). \tag{25}$$

Applying bounded convergence theorem and using (24) and (25), we get that the finite margins of $\{(\widetilde{S}_h(t + \varepsilon u) - \widetilde{S}_h(t))/\varepsilon^{h(t)}\}_{u \geq 0}$ converge to those of $\mathcal{T}(t, \cdot)$ as $\varepsilon \to 0$. It remains to prove the tightness. Let $u \leq v$ such that $|u - v| \leq 1$. We have

$$\mathbb{E}\left[\left(\frac{\widetilde{S}_h(t+\varepsilon u) - \widetilde{S}_h(t)}{\varepsilon^{h(t)}} - \frac{\widetilde{S}_h(t+\varepsilon v) - \widetilde{S}_h(t)}{\varepsilon^{h(t)}}\right)^2\right]$$

$$= \frac{1}{\varepsilon^{2h(t)}}\mathbb{E}[(\widetilde{S}_h(t+\varepsilon u) - \widetilde{S}_h(t+\varepsilon v))^2] = \frac{1}{\varepsilon^{2h(t)}} \int_{t+\varepsilon v}^{t+\varepsilon u} d\theta \int_{t+\varepsilon v}^{t+\varepsilon u} d\sigma \, \mathcal{R}^{*}(\theta, \sigma)$$

$$\leq \sup |R| \int_v^u d\theta \int_v^u d\sigma \, |\theta - \sigma|^{2a-2} \leq \frac{\sup |R|}{2a(2a-1)} |u - v|^{2a}.$$

Then by Kolmogorov criterium, the tightness holds.

Now we assume that $t \neq s$. We have

$$\mathcal{R}^{\#}(t, s, \theta, \sigma, \varepsilon) = \varepsilon^{2-h(t)-h(s)} \mathcal{R}(t + \varepsilon\theta, s + \varepsilon\sigma, h(t + \varepsilon\theta), h(s + \varepsilon\sigma))$$

$$\times |t - s + \varepsilon(\theta - \sigma)|^{h(t+\varepsilon\theta)+h(s+\varepsilon\sigma)-2}. \tag{26}$$

Because $t \neq s$ we can check that for every $\theta$ and $\sigma$ we have as $\varepsilon \to 0$ that

$$\mathcal{R}^{\#}(t, s, \theta, \sigma, \varepsilon) \to 0, \tag{27}$$

and for $\varepsilon$ sufficiently small, $\mathcal{R}^{\#}(t, s, \theta, \sigma, \varepsilon)$ is uniformly bounded for $(\theta, \sigma) \in [0, u] \times [0, v]$. Then by the bounded convergence theorem, we conclude the proof. $\square$

It is classical to deduce pointwise Hölder continuity from local self-similarity [5].

**Proposition 2.** *The process $\widetilde{S}_h$ admits a continuous modification. Moreover, if we assume that $h$ is Hölder continuous, then for every $t_0 \in \mathbb{R}_+$ the pointwise Hölder exponent $\alpha_{\widetilde{S}_h}(t_0)$ of $\widetilde{S}_h$ is almost surely equal to $h(t_0)$.*

**Proof.** We deduce as in [5] from Proposition 1 that $\alpha_{\widetilde{S}_h}(t_0) \leq h(t_0)$. Now we prove that $\alpha_{\widetilde{S}_h}(t_0) \geq h(t_0)$. We let $0 < \eta \leq 1/2$. For every $s$ and $t \in [t_0 - \eta, t_0 + \eta]$ such that $s < t$, we have

$$\mathbb{E}[(\widetilde{S}_h(t) - \widetilde{S}_h(s))^2] \leq \sup |R| \int_s^t d\theta \int_s^t d\sigma \, |\theta - \sigma|^{2\inf_{[s,t]} h - 2}$$

$$\leq \sup |R| |t - s|^{2\inf_{[s,t]} h} \leq \sup |R| |t - s|^{2\inf_{[t_0-\eta, t_0+\eta]} h}.$$

By the fact that $\widetilde{S}_h$ is Gaussian and applying Kolmogorov theorem [7], we get that $\alpha_{\widetilde{S}_h}(t_0) \geq \inf_{[t_0-\eta,t_0+\eta]} h$ for every $0 < \eta \leq 1/2$. Then letting $\eta \to 0$ and using the continuity of $h$, we get $\alpha_{\widetilde{S}_h}(t_0) \geq h(t_0)$, and hence $\alpha_{\widetilde{S}_h}(t_0) = h(t_0)$. □

Please remark these properties are true even if (4) is not fulfilled.

## 5. A representation of the limit process

The aim of this subsection is to give a representation of the limit process $\widetilde{S}_h$ obtained from Theorem 2. This representation uses a universal Gaussian process that we introduce in the following section.

### 5.1. An universal Gaussian field

Here we consider a centered Gaussian field $X = \{X_n(H)\}_{n,H}$. We define $S^N(t,H)$ for every $N$, $t$ and $H$ by

$$S^N(t,H) = \frac{1}{N^H} \sum_{n=1}^{\lfloor Nt \rfloor} X_n(H). \tag{28}$$

**Theorem 3.** *Let $R$ be an asymptotic covariance and $\{X_n(H)\}_{n,H} \in \mathcal{G}_R$ be a centered Gaussian field. Then as $N$ goes to $\infty$, the finite-dimensional distributions of $\{S^N(t,H)\}_{t,H}$ converge to those of a centered Gaussian field $\{\widetilde{W}(t,H)\}_{t,H}$ with covariance given for every $H_1$, $H_2$, $t$ and $s$ by*

$$\begin{aligned}
\mathbb{E}[\widetilde{W}(t,H_1)\widetilde{W}(s,H_2)] &= \frac{R(H_2,H_1)}{(H_1+H_2)(H_1+H_2-1)} s^{H_1+H_2} \\
&\quad + \frac{R(H_1,H_2)}{(H_1+H_2)(H_1+H_2-1)} t^{H_1+H_2} \\
&\quad - \frac{\mathcal{R}(t,s;H_1,H_2)}{(H_1+H_2)(H_1+H_2-1)} |t-s|^{H_1+H_2}.
\end{aligned} \tag{29}$$

**Proof.** We have for every $H_1, H_2, t, s$,

$$\mathbb{E}[S^N(t,H_1)S^N(s,H_2)] = \frac{1}{N^{H_1+H_2}} \sum_{j=1}^{\lfloor Nt \rfloor} \sum_{k=1}^{\lfloor Ns \rfloor} \mathbb{E}[X_k(H_1)X_j(H_2)],$$

and we let

$$K(s,t) = \int_0^t d\theta \int_0^s d\sigma\, \mathcal{R}(\theta,\sigma;H_1,H_2) |\theta-\sigma|^{H_1+H_2-2}.$$

Let $\eta > 0$. Using (12) and $a \leq h \leq b$, there exists an integer $M = M_\eta$ such that for $|j-k| > M$ and every $N$,

$$\frac{1}{N^2} \mathcal{R}(j,k;H_1,H_2) \left|\frac{j-k}{N}\right|^{H_1+H_2-2} - \frac{\eta}{N^2}\left|\frac{j-k}{N}\right|^{H_1+H_2-2}$$

$$\leq \mathbb{E}[X_k(H_1)X_j(H_2)]$$

$$\leq \frac{1}{N^2}\mathcal{R}(j,k;H_1,H_2)\left|\frac{j-k}{N}\right|^{H_1+H_2-2} + \frac{\eta}{N^2}\left|\frac{j-k}{N}\right|^{H_1+H_2-2}. \tag{30}$$

Because of Lemma 2,

$$K(s,t) - \eta L(s,t) \leq \liminf_{N \to \infty} \sum_{j=1}^{\lfloor Nt \rfloor} \sum_{k=1}^{\lfloor Ns \rfloor} \mathbb{E}[X_k(H_1)X_j(H_2)]1_{|j-k|>M}$$

$$\leq \limsup_{N \to \infty} \sum_{j=1}^{\lfloor Nt \rfloor} \sum_{k=1}^{\lfloor Ns \rfloor} \mathbb{E}[X_k(H_1)X_j(H_2)]1_{|j-k|>M} \leq K(s,t) + \eta L(s,t), \quad (31)$$

where

$$L(s,t) = \int_0^t d\theta \int_0^s d\sigma |\theta - \sigma|^{H_1+H_2-2}. \quad (32)$$

We can conclude as the end of the proof of Theorem 2, we get that

$$\lim_{N \to \infty} \mathbb{E}[S^N(t,H_1)S^N(s,H_2)] = \int_0^t d\theta \int_0^s d\sigma \, |\theta - \sigma|^{H_1+H_2-2} \mathcal{R}(\theta,\sigma; H_1, H_2).$$

Then by a direct computation of these last integrals, we get the convergence to a covariance given by (29). □

It is classical to give an alternative form of Theorem 1 which is based on a renormalization group (see, for instance, [15], pages 338–339). Now we propose this alternative approach for Theorem 3.

We fix an asymptotic covariance $R$. Because of Theorem 3, there exists a process $\widetilde{W} = \{\widetilde{W}(t,H)\}_{t,H}$, which is unique in distribution, such that for every $\{X_n(H)\}_{n,H} \in \mathcal{G}_R$

$$\lim_{N \to \infty}^{\mathcal{D}} \left\{ \frac{1}{N^H} \sum_{n=1}^{\lfloor Nt \rfloor} X_n(H) \right\}_{t,H} = \{\widetilde{W}(t,H)\}_{t,H}.$$

We have the equality in distribution for every $\alpha > 0$

$$\{\widetilde{W}(\alpha t, H)\}_{t,H} = \{\alpha^H \widetilde{W}(t,H)\}_{t,H}. \quad (33)$$

We define the field $Z = \{Z_n(H)\}_{n,H}$ for every $n$ and $H$ by

$$Z_n(H) = \widetilde{W}(n,H) - \widetilde{W}(n-1,H).$$

It can be verified that $\{Z_n(H)\}_{n,H} \in \mathcal{G}_R$.

Now we define the renormalization semi-group $T_N$ for every $N$. We let for every $X \in \mathcal{G}_R$

$$T_N X = \{(T_N X)_n(H)\}_{n,H}, \quad (34)$$

where for every $n, H$,

$$(T_N X)_n(H) = \frac{1}{N^H} \sum_{j=nN+1}^{(n+1)N} X_n(H). \quad (35)$$

Because of (33), $Z$ is a fixed point in $\mathcal{G}_R$ of $T_N$ for every $N$. Moreover, for every $X$ in $\mathcal{G}_R$ we have by Theorem 3

$$\lim_{N \to \infty}^{\mathcal{D}} \{(T_N X)_n(H)\}_{n,H} = \{Z_n(H)\}_{n,H}.$$

Finally, we have proved the following result.

**Theorem 4.** *The renormalization semi-group $T_N$ admits the field $Z$ as unique fixed point in the space $\mathcal{G}_R$.*

### 5.2. Representation theorem

Here we consider an asymptotic covariance $R$, a field $X \in \mathcal{G}_R$ and a continuous multifractional function $h : \mathbb{R}_+ \to [a, b] \subset (1/2, 1)$. Recall that by Theorem 2, there exists a multifractional process $\widetilde{S}_h$ such that

$$\lim_{N \to \infty}^{\mathcal{D}} \left\{ \sum_{n=1}^{\lfloor Nt \rfloor} \frac{X_n(h_n^N)}{N^{h_n^N}} \right\}_t = \{\widetilde{S}_h(t)\}_t,$$

and thanks to Theorem 3, there exists a field $\widetilde{W}$ such that

$$\lim_{N \to \infty}^{\mathcal{D}} \left\{ \frac{1}{N^H} \sum_{n=1}^{\lfloor Nt \rfloor} X_n(H) \right\}_{t,H} = \{\widetilde{W}(t, H)\}_{t,H}.$$

In this section, a representation theorem establishes the link between the process $\widetilde{S}_h$ and the field $\widetilde{W}$.

**Theorem 5.** *We assume that the function $R$ is three times continuously differentiable and the function $h$ is two times continuously differentiable. Then we have the following equality in distribution*

$$\{\widetilde{S}_h(t)\}_{t \geq 0} = \left\{ \widetilde{W}(t, h(t)) - \int_0^t h'(\theta) \frac{\partial \widetilde{W}(\theta, H)}{\partial H} \bigg|_{H = h(\theta)} \mathrm{d}\theta \right\}_{t \geq 0}, \tag{36}$$

*where the right-hand side of* (36) *is always defined.*

**Proof.** We can deduce from the assumptions of the theorem and Kolmogorov's criterium that the sample paths of $(t, H) \mapsto \widetilde{W}(t, H)$ are (almost surely and up to a modification) continuous with respect to the first variable and two times continuously differentiable with respect to the second variable. Hence, in particular, the right-hand side of (36) is always defined.

Now we deal with (36). We consider the field $\{Z_n(H)\}_{n,H}$ defined in previous subsection for every $n$ and $H$ by $Z_n(H) = \widetilde{W}(n, H) - \widetilde{W}(n - 1, H)$. On one hand, because of Theorem 2, we have

$$\lim_{N \to \infty}^{\mathcal{D}} \left\{ \sum_{n=1}^{\lfloor Nt \rfloor} \frac{Z_n(h_n^N)}{N^{h_n^N}} \right\}_{t \geq 0} = \{\widetilde{S}_h(t)\}_{t \geq 0}. \tag{37}$$

On the other hand, we shall prove that

$$\lim_{N \to \infty}^{\mathcal{D}} \left\{ \sum_{n=1}^{\lfloor Nt \rfloor} \frac{Z_n(h_n^N)}{N^{h_n^N}} \right\}_{t \geq 0} = \left\{ \widetilde{W}(t, h(t)) - \int_0^t h'(\theta) \frac{\partial \widetilde{W}}{\partial H}(\theta, h(\theta)) \, \mathrm{d}\theta \right\}_{t \geq 0}, \tag{38}$$

which, combined with (37), proves (36). We have

$$\left\{ \sum_{n=1}^{\lfloor Nt \rfloor} \frac{Z_n(h_n^N)}{N^{h_n^N}} \right\}_{t \geq 0} = \left\{ \sum_{n=1}^{\lfloor Nt \rfloor} \frac{1}{N^{h_n^N}} (\widetilde{W}(n, h_n^N) - \widetilde{W}(n - 1, h_n^N)) \right\}_{t \geq 0}$$

$$\stackrel{\mathcal{D}}{=} \left\{ \sum_{n=1}^{\lfloor Nt \rfloor} \left( \widetilde{W}\left(\frac{n}{N}, h_n^N\right) - \widetilde{W}\left(\frac{n - 1}{N}, h_n^N\right) \right) \right\}_{t \geq 0}$$

$$= \{I_1^N(t) - I_2^N(t)\}_{t \geq 0}, \tag{39}$$

where for every $t$ and $N$

$$I_1^N(t) = \sum_{n=1}^{\lfloor Nt \rfloor} \left( \widetilde{W}\left(\frac{n}{N}, h_n^N\right) - \widetilde{W}\left(\frac{n-1}{N}, h_{n-1}^N\right) \right) = \widetilde{W}\left(\frac{\lfloor Nt \rfloor}{N}, h_{\lfloor Nt \rfloor}^N\right),$$

and

$$I_2^N(t) = \sum_{n=1}^{\lfloor Nt \rfloor} \left( \widetilde{W}\left(\frac{n-1}{N}, h_n^N\right) - \widetilde{W}\left(\frac{n-1}{N}, h_{n-1}^N\right) \right).$$

Almost surely,

$$\lim_{N \to \infty} I_1^N(t) = \widetilde{W}(t, h(t))$$

and thanks to Lemma 3 and regularity of $\widetilde{W}$ and $h$,

$$\lim_{N \to \infty} I_2^N(t) = \int_0^t h'(\theta) \frac{\partial \widetilde{W}}{\partial H}(\theta, h(\theta)) \, d\theta.$$

Combined with (39), this proves (38) and concludes the proof. □

## 6. Examples

In this section, we give some examples of multifractional processes $\widetilde{S}_h$ that we can obtain from Theorem 2 as limits of sequences $\{S_h^N\}_N$ for a multifractional function $h \colon [0, \infty) \to [a, b] \subset (1/2, 1)$. In all this section, we assume that $h$ is Hölder continuous.

### 6.1. Fractional white noise model

Let us first consider the case of the multifractional Brownian motion, which at first motivated this article. In constrast to fractional Gaussian noise, which is the fixed point of the renormalization semi-group, it is not the case for the increments of multifractional Brownian motion. In this example, we investigate what is the limit for increments of multifractional Brownian motion. More precisely, let us consider for every $t \in \mathbb{R}$

$$W_H(t) = \frac{1}{C(H)} \int_{-\infty}^{\infty} \frac{e^{-itx} - 1}{|x|^{H+1/2}} \widehat{B}(dx), \tag{40}$$

where $\widehat{B}$ is the Fourier transform of a real Gaussian measure $B$ and the constant $C(H)$ can be written as

$$C(H)^2 = \int_{-\infty}^{\infty} \frac{|e^{-ix} - 1|^2}{|x|^{2H+1}} \, dx = \frac{\pi}{H\Gamma(2H)\sin(\pi H)}.$$

Please note that for each $H$, $W_H$ is a standard fractional Brownian motion. Moreover, if $h \colon [0, \infty) \to (1/2, 1)$ is a multifractional function, then $t \mapsto W_{h(t)}(t)$ is a multifractional Brownian motion. We let

$$X_n(H) = W_H(n+1) - W_H(n).$$

We compute the covariance between $X_j(H_1)$ and $X_k(H_2)$ for every $j$, $k$, $H_1$ and $H_2$:

$$\mathbb{E}[X_j(H_1)X_k(H_2)] = \frac{C((H_1+H_2)/2)^2}{C(H_1)C(H_2)} |j-k|^{H_1+H_2} \tag{41}$$

$$\times \left( \left|1 + \frac{1}{j-k}\right|^{H_1+H_2} + \left|1 - \frac{1}{j-k}\right|^{H_1+H_2} - 2 \right). \tag{42}$$

By a Taylor formula, we get as $u \to 0$

$$|1+u|^{H_1+H_2} = 1 + (H_1+H_2)u + \frac{1}{2}(H_1+H_2)(H_1+H_2-1)u^2 + \mathcal{O}(u^3), \tag{43}$$

where the $\mathcal{O}$ is uniform in $(H_1, H_2) \in [a,b]^2$ because $1/2 < a < b < 1$. Combining (41) and (43), we get that the asymptotic covariance $R$ of $\{X_n(H)\}_{n,H}$ can be written as

$$R(H_1, H_2) = (H_1+H_2)(H_1+H_2-1)\frac{C((H_1+H_2)/2)^2}{C(H_1)C(H_2)}.$$

Applying Theorem 2, we get that $(S_h^N)_N$ converges to the process $\widetilde{S}_h$ with covariance

$$\mathbb{E}[\widetilde{S}_h(t)\widetilde{S}_h(s)] = \int_0^t \mathrm{d}\theta \int_0^s \mathrm{d}\sigma\, |\theta-\sigma|^{h(\theta)+h(\sigma)-2}(h(\theta)+h(\sigma))$$
$$\times (h(\theta)+h(\sigma)-1)\frac{C((h(\theta)+h(\sigma))/2)^2}{C(h(\theta))C(h(\sigma))}.$$

Now we assume that $h$ is continuously differentiable. Using Theorem 5, we can write $\widetilde{S}_h$ as the sum of a multifractional Brownian motion (defined by (2)) and a continuously differentiable process:

$$\widetilde{S}_h(t) = W_{h(t)}(t) - \int_0^t h'(\theta)\frac{\partial W_H(\theta)}{\partial H}\bigg|_{H=h(\theta)} \mathrm{d}\theta. \tag{44}$$

Moreover, we can obtain from (40) the harmonizable representation of the limit process $\widetilde{S}_h$:

$$\widetilde{S}_h(t) = \int_{-\infty}^\infty \bigg\{ \frac{(\mathrm{e}^{\mathrm{i}tx}-1)}{C(h(t))|x|^{h(t)+1/2}}$$
$$- \int_0^t \frac{(\mathrm{e}^{\mathrm{i}\theta x}-1)}{|x|^{h(\theta)+1/2}}\bigg(\frac{\log|x|}{C(h(\theta))} - \frac{C'(h(\theta))}{C(h(\theta))^2}\bigg)h'(\theta)\,\mathrm{d}\theta\bigg\}\widehat{B}(\mathrm{d}x). \tag{45}$$

Note that to establish (45) rigorously, we use the fact that the map $f \mapsto \int_{-\infty}^\infty f(x)\widehat{B}(\mathrm{d}x)$ is an isometry (in particular, $\int_{-\infty}^\infty |f(x)|^2\,\mathrm{d}x = \mathbb{E}[|\int_{-\infty}^\infty f(x)\widehat{B}(\mathrm{d}x)|^2]$ ), and the expression of the covariance of the derivative $\partial W_H/\partial H$ (see, e.g., [1]).

### 6.2. Fractional ARIMA model

In the study of the previous example, we have obtained a symmetric asymptotic covariance $R$. In this section, we present a model for which this symmetric property is not satisfied. We also aim to generalize results in [14], where an invariance principle is established for a class of nonstationary processes with long memory. These processes are generalizations of FARIMA process. At the end of Section 1, the authors wonder if one can extend their results in a continuous time setup. Here we address this question in the Gaussian case.

Let $d \in (0, 1/2)$. We consider the (Gaussian) fractional ARIMA model (FARIMA in short) defined for every $n \in \mathbb{N}$ by (see, for instance, [15])

$$\Phi_n^d = \sum_{l=0}^\infty \frac{\Gamma(d+l)}{l!\Gamma(d)}g_{n-l}, \tag{46}$$

where $\{g_l\}_{l\in\mathbb{Z}}$ are i.i.d. standard normal random variables. Here we let $X_n(H) = \Phi_n^{H-1/2}$. We first establish a lemma.

**Lemma 1.** *As $n \to \infty$, we have*

$$\sup_{d_1,d_2 \in [a-1/2, b-1/2]} \left| n^{1-d_1-d_2} \sum_{l=1}^{\infty} \frac{\Gamma(d_1+n+l)}{(n+l)!} \frac{\Gamma(d_2+l)}{l!} - \int_0^{\infty} (1+u)^{d_1-1} u^{d_2-1} \, \mathrm{d}u \right| \to 0.$$

**Proof.** We have

$$\left| n^{1-d_1-d_2} \sum_{l=1}^{\infty} \frac{\Gamma(d_1+n+l)}{(n+l)!} \frac{\Gamma(d_2+l)}{l!} - \int_0^{\infty} (1+u)^{d_1-1} u^{d_2-1} \, \mathrm{d}u \right|$$
$$\leq \gamma_1(d_1,d_2,n) + \gamma_2(d_1,d_2,n) + \gamma_3(d_1,d_2,n),$$

where

$$\gamma_1(d_1,d_2,n) = \left| n^{1-d_1-d_2} \sum_{l=1}^{\infty} \frac{\Gamma(d_1+n+l)}{(n+l)!} \left\{ \frac{\Gamma(d_2+l)}{l!} - l^{d_2-1} \right\} \right|,$$

$$\gamma_2(d_1,d_2,n) = \left| n^{1-d_1-d_2} \sum_{l=1}^{\infty} l^{d_2-1} \left\{ \frac{\Gamma(d_1+n+l)}{(n+l)!} - (l+n)^{d_1-1} \right\} \right|,$$

$$\gamma_3(d_1,d_2,n) = \left| \frac{1}{n} \sum_{l=1}^{\infty} \left(1+\frac{l}{n}\right)^{d_1-1} \left(\frac{l}{n}\right)^{d_2-1} - \int_0^{\infty} (1+u)^{d_1-1} u^{d_2-1} \, \mathrm{d}u \right|.$$

Let us begin with $\gamma_3$. The function $f_{d_1,d_2} : u \mapsto (1+u)^{d_1-1} u^{d_2-1}$ is decreasing on $[0,\infty)$ so we have

$$\int_{1/n}^{\infty} f_{d_1,d_2}(u) \, \mathrm{d}u \leq \frac{1}{n} \sum_{l=1}^{\infty} f_{d_1,d_2}\left(\frac{l}{n}\right) \leq \int_0^{\infty} f_{d_1,d_2}(u) \, \mathrm{d}u.$$

Hence,

$$\sup_{d_1,d_2} \gamma_3(d_1,d_2,n) \leq \frac{1}{(a-1/2)n^{a-1/2}} \xrightarrow[n \to \infty]{} 0.$$

Now we deal with $\gamma_2$. Using Stirling formula, $\Gamma(z) \sim_{z \to \infty} \sqrt{2\pi} e^{-z} z^{z-1/2}$, we have for every $\eta > 0$ that there exists $M_\eta$ depending only on $\eta$, $a$ and $b$ such that for every $l \geq M_\eta$

$$\left| \frac{\Gamma(d_1+l)}{l!} - l^{d_1-1} \right| \leq \eta l^{d_1-1} e^{1-d_1}.$$

We get for every $n \geq M_\eta$

$$\gamma_2(d_1,d_2,n) \leq n^{1-d_1-d_2} \sum_{l=1}^{\infty} \left| l^{d_2-1} \left\{ \frac{\Gamma(d_1+n+l)}{(n+l)!} - (l+n)^{d_1-1} \right\} \right|$$
$$\leq \eta n^{1-d_1-d_2} \sum_{l=1}^{\infty} l^{d_2-1} (l+n)^{d_1-1}.$$

So, for every $\eta > 0$, $\limsup_{n \to \infty} \sup_{d_1,d_2} \gamma_2(d_1,d_2,n) \leq \eta S$ where $S = \sup_{d_1,d_2} \int_0^{\infty} f_{d_1,d_2}(u) \, \mathrm{d}u < \infty$. Hence, $\lim_{n \to \infty} \sup_{d_1,d_2} \gamma_2(d_1,d_2,n) = 0$. We use a similar argument for $\gamma_1$. □

Using Lemma 1, the function $R$ is

$$R(H_1,H_2) = \Gamma\left(H_1 - \frac{1}{2}\right)^{-1} \Gamma\left(H_2 - \frac{1}{2}\right)^{-1} \int_0^{\infty} (1+u)^{H_1-3/2} u^{H_2-3/2} \, \mathrm{d}u$$

$$= \frac{1}{\pi}\sin\left(\pi\left(H_1-\frac{1}{2}\right)\right)\Gamma(2-H_1-H_2),$$

and the limit process $\widetilde{S}_h$ has the covariance

$$\mathbb{E}[\widetilde{S}_h(t)\widetilde{S}_h(s)] = \int_0^t d\theta \int_0^s d\sigma \frac{|\theta-\sigma|^{h(\theta)+h(\sigma)-2}}{\Gamma(h(\theta)-1/2)\Gamma(h(\sigma)-1/2)}$$
$$\times \int_0^\infty (1+u)^{h(\theta\vee\sigma)-3/2} u^{h(\theta\wedge\sigma)-3/2}\,du$$
$$= \frac{1}{\pi} \int_0^t d\theta \int_0^s d\sigma\, |\theta-\sigma|^{h(\theta)+h(\sigma)-2} \sin\left(\pi\left(h(\theta\vee\sigma)-\frac{1}{2}\right)\right)$$
$$\times \frac{\Gamma(2-h(\theta)-h(\sigma))}{\Gamma(h(\theta)-1/2)\Gamma(h(\sigma)-1/2)}.$$

**Appendix. Riemann sum convergence type lemma**

**Lemma 2.** *Let $t, s \geq 0$, $G: [0,t] \times [0,s] \to \mathbb{R}$ be a continuous function and $M > 1$. We have*

$$\lim_{N\to\infty} \frac{1}{N^2} \sum_{j=1}^{\lfloor Nt\rfloor} \sum_{k=1}^{\lfloor Ns\rfloor} G\left(\frac{j}{N},\frac{k}{N}\right) \left|\frac{j-k}{N}\right|^{h_k^N+h_j^N-2} 1_{|j-k|>M} = \int_0^t d\theta \int_0^s d\sigma\, G(\theta,\sigma)|\theta-\sigma|^{h(\theta)+h(\sigma)-2}, \quad (47)$$

*where the integral in the right-hand side of (47) is always defined.*

**Proof.** We let for every $j,k,N$ and $(\theta,\sigma) \in (\frac{j-1}{N},\frac{j}{N}] \times (\frac{k-1}{N},\frac{k}{N}]$

$$G_N(\theta,\sigma) := G\left(\frac{j}{N},\frac{k}{N}\right) \left|\frac{j-k}{N}\right|^{h_k^N+h_j^N-2} 1_{|j-k|>M}.$$

Notice that

$$\int_0^t d\theta \int_0^s d\sigma\, G_N(\theta,\sigma) = \frac{1}{N^2} \sum_{j=1}^{\lfloor Nt\rfloor} \sum_{k=1}^{\lfloor Ns\rfloor} G\left(\frac{j}{N},\frac{k}{N}\right) \left|\frac{j-k}{N}\right|^{h_k^N+h_j^N-2} 1_{|j-k|>M},$$

hence, it suffices to prove

$$\lim_{N\to\infty} \int_0^t d\theta \int_0^s d\sigma\, G_N(\theta,\sigma) = \int_0^t d\theta \int_0^s d\sigma\, G(\theta,\sigma)|\theta-\sigma|^{h(\theta)+h(\sigma)-2}. \quad (48)$$

We shall use a dominated convergence argument. Since $G$ is continuous, then for almost every $(\theta,\sigma) \in [0,t] \times [0,s]$,

$$\lim_{N\to\infty} G_N(\theta,\sigma) = G(\theta,\sigma)|\theta-\sigma|^{h(\theta)+h(\sigma)-2}. \quad (49)$$

It remains to establish the boundedness of $|G_N|$ by an integrable function. Since $G$ is bounded, then for every $(\theta,\sigma) \in (\frac{j-1}{N},\frac{j}{N}] \times (\frac{k-1}{N},\frac{k}{N}]$

$$|G_N(\theta,\sigma)| \leq \sup|G|\left(\left|\frac{j-k}{N}\right|^{2a-2} + \left|\frac{j-k}{N}\right|^{2b-2}\right)$$

and since $|j-k| > M > 1$,

$$|\theta - \sigma| \leq \left|\frac{j-k}{N}\right| + \frac{1}{N} \leq 2\left|\frac{j-k}{N}\right|,$$

hence, for almost every $(\theta, \sigma) \in [0,t] \times [0,s]$

$$|G_N(\theta, \sigma)| \leq \sup|G|(|\theta - \sigma|^{2a-2} + |\theta - \sigma|^{2b-2}).$$

The map $(\theta, \sigma) \mapsto |\theta - \sigma|^{2a-2} + |\theta - \sigma|^{2b-2}$ is integrable on $[0,t] \times [0,s]$, then thanks to (49) and dominated convergence, we get (48), which proves (47). $\square$

**Lemma 3.** *Let $a < b \in \mathbb{R}$ and two functions $f \colon \mathbb{R}_+ \times [a,b] \to \mathbb{R}$ and $h \colon \mathbb{R}_+ \to [a,b]$. We assume that $h$ is two times continuously differentiable and $f$ is two times continuously differentiable with respect to its second variable. Then for every $t \geq 0$, we have*

$$\lim_{N \to \infty} \sum_{n=1}^{\lfloor Nt \rfloor} \left( f\left(\frac{n-1}{N}, h\left(\frac{n}{N}\right)\right) - f\left(\frac{n-1}{N}, h\left(\frac{n-1}{N}\right)\right) \right) = \int_0^t h'(\theta)\, \partial_2 f(\theta, h(\theta))\, \mathrm{d}\theta.$$

**Proof.** By the Taylor formula, for every $n$ and $N$, there exists $t_{n,N} \in [(n-1)/N, n/N]$ such that

$$f\left(\frac{n-1}{N}, h\left(\frac{n}{N}\right)\right) - f\left(\frac{n-1}{N}, h\left(\frac{n-1}{N}\right)\right)$$
$$= \frac{1}{N} h'\left(\frac{n-1}{N}\right) \partial_2 f\left(\frac{n-1}{N}, h\left(\frac{n-1}{N}\right)\right) + \frac{1}{N^2} E_{n,N}, \tag{50}$$

where

$$E_{n,N} = h''(t_{n,N})\, \partial_2 f\left(\frac{n-1}{N}, h(t_{n,N})\right) + (h'(t_{n,N}))^2\, \partial_2^2 f\left(\frac{n-1}{N}, h(t_{n,N})\right).$$

Thanks to the regularity properties of $h$ and $f$, there exists $C > 0$ such that for every $n$ and $N$, $|E_{n,N}| \leq C$. Then using (50),

$$\left| \sum_{n=1}^{\lfloor Nt \rfloor} \left( f\left(\frac{n-1}{N}, h\left(\frac{n}{N}\right)\right) - f\left(\frac{n-1}{N}, h\left(\frac{n-1}{N}\right)\right) \right) - \frac{1}{N} \sum_{n=1}^{\lfloor Nt \rfloor} h'\left(\frac{n-1}{N}\right) \partial_2 f\left(\frac{n-1}{N}, h\left(\frac{n-1}{N}\right)\right) \right| \leq \frac{Ct}{N}. \tag{51}$$

We conclude the proof by combining classical Riemann sums convergence with (51). $\square$